# Error-rate and decision-theoretic methods of multiple testing

## Which genes have high objective probabilities of differential expression?

**David R. Bickel**

*Office of Biostatistics and Bioinformatics\**
*Medical College of Georgia*
*Augusta, GA 30912-4900*

*bickel@prueba.info*
*www.davidbickel.com*

**Abstract.**
   Given a multiple testing situation, the null hypotheses that appear to have sufficiently low probabilities of truth may be rejected using a simple, nonparametric method of decision theory. This applies not only to posterior levels of belief, but also to conditional probabilities in the sense of relative frequencies, as seen from their equality to local false discovery rates (dFDRs). This approach neither requires the estimation of probability densities, nor of their ratios. Decision theory can inform the selection of false discovery rate weights. Decision theory is applied to gene expression microarrays with discussion of the applicability of the assumption of weak dependence.

\*Address after March 31, 2004: 7250 NW 62nd St., P. O. Box 552, Johnston, IA 50131-0552

# 1 Introduction

In multiple hypothesis testing, the ratio of the number of false discoveries (false positives) to the total number of discoveries is often of interest. Here, a *discovery* is the rejection of a null hypothesis and a *false discovery* is the rejection of a true null hypothesis. To quantify this ratio statistically, consider $m$ hypothesis tests with the $i$th test yielding $H_i = 0$ if the null hypothesis is true or $H_i = 1$ if it is false and the corresponding alternative hypothesis is true. Let $R_i = 1$ if the $i$th null hypothesis is rejected and $R_i = 0$ otherwise and define $V_i = (1 - H_i) R_i$, such that $\sum_{i=1}^{m} R_i$ is the total number of discoveries and $\sum_{i=1}^{m} V_i$ is the total number of false discoveries. The intuitively appealing quantity $\mathrm{E}(\sum_{i=1}^{m} V_i / \sum_{i=1}^{m} R_i)$ is not a viable definition of a false discovery rate since it usually cannot be assumed that $\mathrm{P}(\sum_{i=1}^{m} R_i > 0) = 1$. Benjamini and Hochberg (1995) listed three possible candidates for such a definition:

$$\varepsilon \equiv \mathrm{E}\left(\frac{\sum_{i=1}^{m} V_i}{\sum_{i=1}^{m} R_i} \,\bigg|\, \sum_{i=1}^{m} R_i > 0\right) \mathrm{P}\left(\sum_{i=1}^{m} R_i > 0\right); \tag{1}$$

$$Q \equiv \mathrm{E}\left(\frac{\sum_{i=1}^{m} V_i}{\sum_{i=1}^{m} R_i} \,\bigg|\, \sum_{i=1}^{m} R_i > 0\right); \tag{2}$$

$$\nabla \equiv \frac{\mathrm{E}(\sum_{i=1}^{m} V_i)}{\mathrm{E}(\sum_{i=1}^{m} R_i)}. \tag{3}$$

Benjamini and Hochberg (1995) called $\varepsilon$ the *false discovery rate* (FDR) since null hypotheses can be rejected in such a way that $\varepsilon \leq \alpha$ is guaranteed for any significance level $\alpha \in (0, 1)$. This *control* of the FDR is accomplished by choosing a suitable rejection region $\Gamma \subset S$, where $S$ is the union of the test-statistic space and the empty set, and by rejecting the $i$th null hypothesis if and only if the $i$th test statistic, $t_i \in S$, is an element of $\Gamma$. Benjamini and Hochberg (1995) expressed this in terms of p-values instead of test statistics, with $S=[0,1]$. The FDR was originally controlled under independence of the test statistics (Benjamini and Hochberg, 1995), but can also be controlled for more general cases (Benjamini and Yekutieli, 2001). The other two candidates, $Q$ and $\nabla$, cannot be controlled in this way since, if all $m$ of the null hypotheses are true, then $Q = 1 > \alpha$ and $\Delta = 1 > \alpha$ (Benjamini and Hochberg, 2002). Nonetheless, Storey (2002a) recommended the use of $Q$, which he called the *positive false discovery rate* (pFDR), since the multiplicand $P(\sum_{i=1}^{m} R_i > 0)$ of the FDR makes it harder to interpret than the pFDR. Although $\nabla$ of Eq. (3) is also easier to interpret than the FDR, Storey (2003) complained that $\nabla$ fails to describe the simultaneous fluctuations in $\sum_{i=1}^{m} V_i$ and $\sum_{i=1}^{m} R_i$, in spite of its attraction as a simple measure. However, Fernando *et al*. (2004) demonstrated that $\nabla$, the *proportion of false positives* (PFP), is convenient for multiple sets of comparisons, when the FDR and pFDR are less intuitive. When $H_i$ is considered as a random variable, the PFP becomes the "Bayesian FDR" of Efron and Tibshirani (2002), but frequentist inference also applies to the PFP. For example, the PFP is often equal to a frequentist probability that a rejected null hypothesis is true (Fernando *et al*. 2004). A drawback of the PFP is that it is undefined for $E(\sum_{i=1}^{m} R_i) = 0$, which can occur if all null hypotheses are false, and which prevents the PFP from being controlled in the sense that the FDR can be controlled. This disadvantage is overcome by introducing an error rate that is equal to the PFP if the PFP is defined, and equal to 0 otherwise (Bickel, 2004a):

$$\Delta \equiv \begin{cases} \frac{E(\sum_{i=1}^{m} V_i)}{E(\sum_{i=1}^{m} R_i)} = \nabla & \text{if} \quad P(\sum_{i=1}^{m} R_i > 0) \neq 0 \\ 0 & \text{if} \quad P(\sum_{i=1}^{m} R_i > 0) = 0 \end{cases}. \tag{4}$$

Unlike the PFP, this error rate equals the probabity that a rejected null hypothesis is true, even when each null hypothesis is almost never rejected. (A probability conditional on an event of zero probabilty is defined arbitrarily as zero (Breiman 1992).) It will be seen that, unlike the error rates of (1)-(3), $\Delta$ can be optimized using decision theory without dependence assumptions, optimization that rejects sufficiently improbable null hypotheses under general conditions. Because of this advantage in a decision-theoretic framework, $\Delta$ has been called the *decisive false discovery rate* (dFDR; Bickel, 2004a). In addition, the dFDR has desirable properties independent of decision theory: it can be estimated using estimators of the FDR or the PFP, and it can be controlled in the same sense as the FDR, i.e., one may reject as many null hypotheses as possible with the constraint that $\Delta \leq \alpha$, even when none of the null hypotheses are true (Bickel, 2004a).

## 2 Decision theory

### 2.1 Desirability and the dFDR

Denote by $b_i$ the *benefit* of rejecting the *i*th null hypothesis when it is false, where the benefit can be economic, such as an amount of money, or some other positive desirability in the sense of Jeffrey (1983). Let $c_i$ be the *cost* of rejecting the *i*th null hypothesis when it is true, with cost here playing the role of a negative desirability, such that we have $b_i \geq 0$, $c_i \geq 0$, and the vectors $\boldsymbol{b} \equiv (b_i)_{i=1}^m$ and $\boldsymbol{c} \equiv (c_i)_{i=1}^m$. Then the net desirability is

$$d(\boldsymbol{b}, \boldsymbol{c}) = \sum_{i=1}^m (b_i(R_i - V_i) - c_i V_i) = \sum_{i=1}^m b_i R_i - \sum_{i=1}^m (b_i + c_i) V_i. \tag{5}$$

As seen in Section 4.2, this formulation in terms of both costs and benefits gives the same results as the equivalent formulation in terms of a loss function of costs of Type I and Type II errors. If the costs and benefits are independent of the hypothesis, then the net desirability is $d = b_1 \left( \sum_{i=1}^m R_i - \left(1 + \frac{c_1}{b_1}\right) \sum_{i=1}^m V_i \right) = b_1 \left(1 - \left(1 + \frac{c_1}{b_1}\right) \sum_{i=1}^m V_i / \sum_{i=1}^m R_i \right) \sum_{i=1}^m R_i$ and its expectation value is

$$\mathrm{E}(d) \approx b_1 \left(1 - \left(1 + \frac{c_1}{b_1}\right) Q \right) \mathrm{E}\left( \sum_{i=1}^m R_i \right) \tag{6}$$

if $Q \approx \Delta$. (Genovese and Wasserman (2001) instead utilized a marginal risk function that is effectively equivalent to $-\mathrm{E}(d)/m$ if $b_1 = c_1$.) Thus, the rejection region can be chosen such that the pFDR (or FDR) and the expected number of rejections maximize the approximate expected desirability, given a cost-to-benefit ratio. However, these approximations are not needed since there is an exact result in terms of the PFP:

$$E(d) = b_1 E\left(\sum_{i=1}^{m} R_i - \left(1 + \frac{c_1}{b_1}\right) \sum_{i=1}^{m} V_i\right) = \begin{cases} b_1 \left(1 - \left(1 + \frac{c_1}{b_1}\right) \nabla\right) E(\sum_{i=1}^{m} R_i) & \text{if } E(\sum_{i=1}^{m} R_i) \neq 0 \\ 0 & \text{if } E(\sum_{i=1}^{m} R_i) = 0 \end{cases}, \quad (7)$$

more succinctly expressed as a function of the dFDR:

$$E(d) = b_1 \left(1 - \left(1 + \frac{c_1}{b_1}\right) \Delta\right) E\left(\sum_{i=1}^{m} R_i\right). \tag{8}$$

When all of the null hypotheses are true, $\Delta = 1$ unless $E(\sum_{i=1}^{m} R_i) = 0$, in which case $\Delta = 0$. Likewise, when all of the null hypotheses are false, $\Delta = 0$ for all rejection regions and $E(\sum_{i=1}^{m} R_i) = m$ is obviously optimal; in this case, the rejection region would be the whole test statistic space.

## 2.2 Conditional probability that a null hypothesis is true

Whether or not $H_i$ is a random variable, the rejection region that maximizes $E(d)$ is the same rejection region that selects which null hypotheses have probabilities of truth as low as or lower than some value. To make the terms of equation (7) more succinct and explicitly dependent on the rejection region, let $D(\Gamma) = E(d)$, $R_\bullet(\Gamma) = E(\sum_{i=1}^{m} R_i)$, $V_\bullet(\Gamma) = E(\sum_{i=1}^{m} V_i)$, and $p = \left(1 + \frac{c_1}{b_1}\right)^{-1}$, so that the optimal rejection region, $\Gamma_{\text{optimal}}$, is defined by

$$D(\Gamma_{\text{optimal}}) = \max_{\Gamma \subset S} D(\Gamma) = \frac{b_1}{p} \max_{\Gamma \subset S} [p R_\bullet(\Gamma) - V_\bullet(\Gamma)]. \tag{9}$$

In other words, given a cost-to-benefit ratio $c_1 / b_1$, $\Gamma_{\text{optimal}}$ is the rejection region that maximizes the relative expected gain $G(\Gamma)$, defined as $p R_\bullet(\Gamma) - V_\bullet(\Gamma)$. This notation facilitates obtaining the result that if it is desirable to reject hypotheses with test statistics in some region $\Gamma$, then the probability of a false rejection is no greater than $p$:

$$D(\Gamma) \geq 0 \Rightarrow G(\Gamma) \geq 0 \Rightarrow P(H_i = 0 \mid t_i \in \Gamma) = \frac{V_\bullet(\Gamma)}{R_\bullet(\Gamma)} \leq p. \tag{10}$$

This inequality holds not only if $H_i$ is random, as in a Bayesian mixture model, but also if, for each value of $i$, $H_i$ is fixed and $t_i$ is random. In this case, $P(H_i = 0 \mid t_i \in \Gamma)$ is not a level of posterior belief, but is the long-term frequency of true null hypotheses relative to rejected null hypotheses for a fixed rejection region $\Gamma$. (That the PFP is related to a frequentist probability was also noted by Fernando *et al.* (2004).) Equation (10) implies that $p$ is also the upper bound on every local rejection region in $\Gamma_{\text{optimal}}$:

$$\forall_{\Gamma \subset \Gamma_{\text{optimal}}} P(H_i = 0 \mid t_i \in \Gamma) \leq p. \tag{11}$$

This can be proven by contradiction as follows. Assume $\exists_{\Gamma \subset \Gamma_{\text{optimal}}} P(H_i = 0 \mid t_i \in \Gamma) > p$. Then, by equation (10), $D(\Gamma) < 0$. But $D(\Gamma) < 0$ means the expected cost exceeds the expected benefit for $\Gamma$, which cannot be true of any subset of a globally optimal rejection region. It follows that $\Gamma \notin \Gamma_{\text{optimal}}$, which contradicts the premise.

A corollary is that the parameter $p$ is an upper bound on the probability of a false positive given a specific value of the test statistic:

$$\begin{aligned} &\forall_{[\tau, \tau+dt) \subset \Gamma_{\text{optimal}}} P(H_i = 0 \mid t_i \in [\tau, \tau+dt)) \leq p \\ &\therefore \forall_{\tau \in \Gamma_{\text{optimal}}} P(H_i = 0 \mid t_i = \tau) \leq p. \end{aligned} \tag{12}$$

Thus, decision-theoretic optimization rejects all null hypotheses with probabilities less than or equal to $p$. This relation holds for $P(H_i = 0 \mid t_i = \tau)$ as a long-term, relative-frequency probability. An equivalent relationship also holds for subjective posterior probabilities (Müller *et al.* 2004; Bickel 2004c). (As Jeffrey (1983) and earlier references in Scott and Berger (2003) show, there is a long history of applying Bayesian statistics to decision theory, but equation (12) appears to be a new result for frequentist inference.) The practical applicaton of equation (12) is that an investigator unable to provide a cost-to-benefit ratio may instead specify a probability theshold $p$ and then compute the ratio as $\frac{c_1}{b_1} = \frac{1}{p} - 1$. For example, analogous to the usual Type I significance level of 5%, optimizing the desirability with a cost-to-benefit ratio of 19 would limit the probability of a false discovery to the 5% level, provided that the expected desirability is positive for some rejection region. The interpretation is straightforward. Simple dFDR control at the 5% level finds as many discoveries as possible such that for every expected false discovery, there are at least 19 expected true discoveries. That leads to less than 19 expected true discoveries per expected false discovery near the border of the rejection region. Decision-theoretic optimization at a cost-to-benefit ratio of 19 instead maximizes the net benefit minus the net cost, given that every false discovery offsets 19 true discoveries. That keeps the local dFDR under 5%, even at the border of the rejection region. In terms of the rejection subregion $\Gamma_{local}$, the *local dFDR* is defined analogously to the local pFDR of Efron *et al.* (2001):

$$\Delta_{local}[\![\Gamma_{local}]\!] \equiv \begin{matrix} \frac{E(\sum_{i=1}^m V_i(\Gamma_{local}))}{E(\sum_{i=1}^m R_i(\Gamma_{local}))} & P(\sum_{i=1}^m R_i(\Gamma_{local}) > 0) \neq 0 \\ 0 & P(\sum_{i=1}^m R_i(\Gamma_{local}) > 0) = 0 \end{matrix};$$

$$\Gamma_{local} \subset \Gamma_{optimal}; \quad \begin{matrix} t_i \in \Gamma_{local} \Rightarrow V_i(\Gamma_{local}) = V_i, R_i(\Gamma_{local}) = R_i \\ t_i \notin \Gamma_{local} \Rightarrow V_i(\Gamma_{local}) = 0, R_i(\Gamma_{local}) = 0 \end{matrix}.$$

Equations (10) and (11) yield

$$\forall_{\Gamma_{local} \subset \Gamma_{optimal}} \Delta_{local}[\![\Gamma_{local}]\!] \leq p. \tag{14}$$

It follows that maximizing the net desirability is equivalent to keeping the local dFDRs as well as the probabilities below a specified threshold. It will be seen that this does not require probability density or density ratio estimation, as does the method of Efron *et al.* (2001). Although the decision-theoretic approach does require estimation of the cumulative probability in the nonparametric case, such estimation is less noisy than density estimation. The tradeoff is a need to perform a maximization, which might be viewed as approximating a derivative.

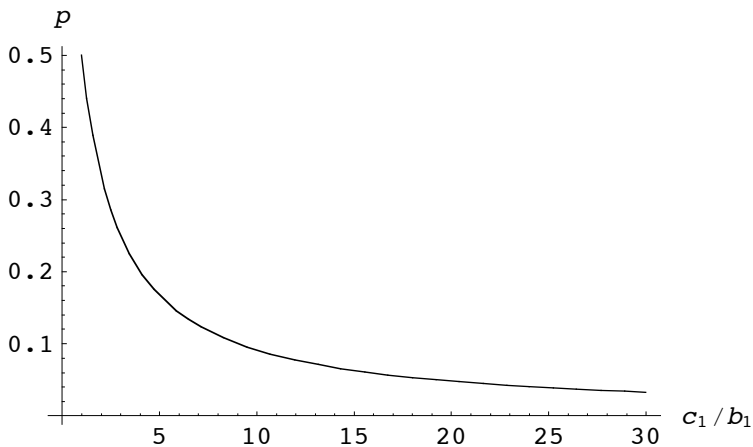

Upper bound of conditional probability for max E($d$).

**Figure 1**

## 3  Estimation of dFDR and optimization of desirability

If $\mathrm{E}(\sum_{i=1}^{m} R_i) \neq 0$, the PFP, and thus the dFDR, can be written in terms of $\pi_0$, the proportion of null hypotheses that are true:

$$\nabla = \frac{\mathrm{E}(\sum_{i=1}^{m} (1 - H_i) R_i)}{\mathrm{E}(\sum_{i=1}^{m} R_i)} = \frac{\pi_0 \sum_{i=1}^{m} \mathrm{P}(t_i \in \Gamma \mid H_i = 0)}{\sum_{i=1}^{m} \mathrm{P}(t_i \in \Gamma)}. \tag{15}$$

Since $\pi_0$ is unknown, it must either be estimated or conservatively set to an upper bound such as 1. The two sums must also be estimated to yield an estimate of $\nabla$. For ease of notation, let $\Gamma = [\tau, \infty)$, and denote by $F$ the distribution of the test statistics and by $F^{(0)}$ the distribution of the test statistics under each null hypothesis. (The same method applies to more general rejection regions, but a single-interval rejection region avoids the bias, mentioned by Storey and Tibshirani (2003), of choosing a multiple-interval region from the data and then using that region to obtain estimates.) Then, if $\mathrm{P}(t \geq \tau) \neq 0$,

$$\nabla(\tau) = \frac{\pi_0 \, m \mathrm{P}(t \geq \tau \mid H = 0)}{m \mathrm{P}(t \geq \tau)} = \frac{\pi_0 (1 - F^{(0)}(\tau))}{1 - F(\tau)}, \tag{16}$$

which, given $B$ bootstrap or permutation samples, is naturally estimated by

$$\hat{\varepsilon}(\tau) \equiv \hat{\nabla}(\tau) \equiv \frac{\hat{\pi}_0 \sum_{i=1}^{mB} I(T_i^{(0)} \geq \tau) \big/ (mB)}{\sum_{i=1}^{m} I(T_i \geq \tau)/m} \quad \text{for} \quad \sum_{i=1}^{m} I(T_i \geq \tau) \neq 0, \tag{17}$$

where $\hat{\pi}_0$ is an estimate of $\pi_0$, $I(x) = 1$ if $x$ is true and $I(x) = 0$ if $x$ is false, $T_i$ is the observed test statistic of the $i$th test, and $T_i^{(0)}$ is the $i$th test statistic out of $M = mB$ test statistics calculated by a resampling procedure under the null hypothesis.

Such resampling is not needed when $F^{(0)}$ is known, as when $-T_i$ is the p-value of the $i$th test; Storey (2002a) and Genovese and Wasserman (2002) examined this special case. Storey (2002a) and Genovese and Wasserman (2002) effectively used additive inverses of p-values as the test statistics, taking advantage of the fact that p-values are uniformly distributed under the null hypothesis. Thus, given that $\mathcal{P}_i$ is the p-value of the $i$th test, equation (20) can be replaced by the simpler dFDR estimator

$$\hat{\varepsilon}_\Pi \equiv \hat{\nabla}_\Pi \equiv \frac{\hat{\pi}_0 \Pi}{\sum_{i=1}^m I(\mathcal{P}_i \leq \Pi)/m} \text{ for } \sum_{i=1}^m I(\mathcal{P}_i \leq \Pi) \neq 0, \quad (18)$$

where $\Pi$ is the p-value rejection region threshold ($R_i = 1$ if $\mathcal{P}_i \leq \Pi$ or $R_i = 0$ if $\mathcal{P}_i > \Pi$). Bickel (2004a) took this approach with Wilcoxon rank-sum p-values.

Making use of the fact that a large proportion of the test statistics that are less than some threshold $\lambda$ would be associated with true null hypotheses, Efron et al. (2001) and Storey (2002a) estimated $\pi_0$ by variants of

$$\hat{\pi}_0(\lambda) = \frac{\sum_{i=1}^m I(T_i < \lambda)/m}{\sum_{i=1}^M I(T_i^{(0)} < \lambda)/M}. \quad (19)$$

Storey (2003) and Storey, Taylor, and Siegmund (2004) offer ways to optimize $\lambda$, but in the present application to gene expression data, the more straightforward approach of Efron et al. (2001) will be followed for simplicity. They transformed the observed test statistics to follow $\Phi$, the standard normal distribution, then counted the numbers of transformed statistics falling in the range $(-1/2, 1/2)$. One could use Eq. (19) to obtain similar results, without transformation, by choosing $\lambda$ such that the proportion of null test statistics that are less than $\lambda$ is as close as possible to $\Phi(1/2) - \Phi(-1/2) \approx 0.382925$. Alternately, one may follow the simpler method of Benjamini and Hochberg (1995): assume that $\pi_0 = 1$ for purposes of estimation. That may be too conservative when a large portion of null hypotheses are false. The estimation of $\pi_0$ does not rely on distributional assumptions, which can require extensive model validation.

Once an estimator of the PFP is chosen, it can be used to estimate the dFDR:

$$\hat{\Delta}(\tau) \equiv \begin{cases} \hat{\nabla}(\tau) & \text{if } \sum_{i=1}^m I(T_i \geq \tau) \neq 0 \\ 0 & \text{if } \sum_{i=1}^m I(T_i \geq \tau) = 0 \end{cases}; \quad (20)$$

$$\hat{\Delta}_\Pi \equiv \begin{matrix} \hat{V}_\Pi & \text{if} & \sum_{i=1}^m I(\mathcal{P}_i \leq \Pi) \neq 0 \\ 0 & \text{if} & \sum_{i=1}^m I(\mathcal{P}_i \leq \Pi) = 0 \end{matrix}. \qquad (21)$$

Equations (20) and (21) give the estimator for statistics of one-sided tests and for p-values, respectively. Like estimation of the FDR and pFDR, estimation of the dFDR does not require restrictive assumptions about $F$. Given a value of $\lambda$ and a test-independent cost-to-benefit ratio, Eq. (8) suggests optimizing the dFDR by selecting the value of $\tau \in \{T_i\}_{i=1}^m$ that maximizes

$$\hat{D}(\tau) = b_1 \left(1 - \left(1 + \frac{c_1}{b_1}\right)\hat{\Delta}(\tau)\right) \sum_{i=1}^m I(T_i \geq \tau), \qquad (22)$$

with $\hat{\Delta}_\Pi$ substituted for $\hat{\Delta}(\tau)$, if desired. Clearly, the best rejection region can be found without knowledge of $b_1$, as long as the ratio $c_1/b_1$ is specified. $\hat{D}(\tau)$ is a conservatively inconsistent estimator of $E(d)$ for any rejection region, given the weak dependence of test statistics or p-values and the other conditions of Storey, Taylor, and Siegmund (2004) for which conservative inconsistency of $\hat{V}(\tau)$ holds. Under those conditions, just as keeping $\hat{\Delta}(\tau)$ less than or equal to some $\alpha$ controls the FDR at level $\alpha$ (Storey, Taylor, and Siegmund 2004), it also controls the dFDR at level $\alpha$.

The R and S-PLUS functions used to perform estimation are available through www.davidbickel.com.

# 4 Some theoretical relationships

## 4.1 dFDR as a Bayesian or frequentist probability

While the above definition of dFDR and estimation procedure were given in terms agreeable with a frequentist perspective in which the truth or falsehood of each null hypothesis is fixed, they are equally valid in a Bayesian construction. Efron et al. (2001), Genovese and Wasserman (2002), and Storey (2002a) considered the case in which each $H_i$ is a random variable and $\pi_0$ is the prior probability that a null hypothesis is true, such that the marginal distribution $F$ is a mixture of $F^{(0)}$ and $F^{(1)}$, the distributions of the statistics conditioned on the truth and falsehood of the null hypothesis: $F = \pi_0 F^{(0)} + (1 - \pi_0) F^{(1)}$. Storey (2003) pointed out that the PFP is equal to the posterior probability that a null hypothesis is true:

$$\nabla = \frac{P(H = 0) P(t \in \Gamma \mid H = 0)}{P(t \in \Gamma)} = \frac{P(H = 0, t \in \Gamma)}{P(t \in \Gamma)} = P(H = 0 \mid t \in \Gamma) \text{ for } P(t \in \Gamma) \neq 0. \quad (23)$$

He also proved that this holds for the pFDR, assuming that the elements of $\{t_i\}_{i=1}^m$ and $\{H_i\}_{i=1}^m$ are independent, with the corollary that $\Delta = Q$ for the same rejection region $\Gamma$. Thus, under approximate independence,

$$Q \approx P(H = 0 \mid t \in \Gamma), \quad (24)$$

with equality only holding under exact independence. Storey (2002a) used relation (24) to motivate the methods of estimating the FDR and pFDR that inspired the variant method of Eqs. (20) and (19). (Storey (2002a) adjusted the estimates for the conditionality of the pFDR on $\sum_{i=1}^m R_i > 0$, but the adjustment is not needed to estimate the FDR, PFP, or dFDR.) Storey (2003) also proved that such estimation procedures are conservatively inconsistent ($m \to \infty$) under weak dependence.

The dFDR is not only equal to a posterior probability, but is equal to a conditional probability without any need for a Bayesian approach. In addition to making the dFDR easy to interpret, this equality, like equation (12), enables the investigator to use of rules of probability to make additional inferences. For example, the rules of adding disjoint probabilities and multiplying independent probabilities aid in drawing conclusions about gene networks studied with the dFDR (Bickel, 2004b).

## 4.2 Other methods of optimization

The maximization of $\hat{D}(\tau)$ determines a rejection region using the dFDR and the cost and benefit of a false or true discovery. Genovese and Wasserman (2001), Storey (2003), Müller *et al.* (2004), and Scott and Berger (2003), on the other hand, presented cost functions that use the concept of a *false nondiscovery*, a failure to reject a false null hypothesis. The approach taken herein has much more in common with the approaches of Storey (2003) and Müller *et al.* (2004) than with that of Genovese and Wasserman (2001), also considered by Storey (2003) and Müller *et al.* (2004). Scott and Berger (2003) presented a parametric Bayesian technique. Each will be described in turn.

Storey (2003) introduced a method based on $c'$, the cost of a false nondiscovery, instead of the benefit from a true discovery. Table 1 clarifies the difference between his formulation and that described above:

| Change in desirability $d_i$ (c, b, c′) | $R_i$ = 0 | $R_i$ = 1 |
|---|---|---|
| $H_i$ = 0 | 0 | $-c$ |
| $H_i$ = 1 | $-c'$ | $+b$ |

Contribution of each cost or benefit to the net desirability.

**Table 1**

The approach of Eqs. (5), (6), and (8) requires that $c > 0$, $b > 0$, and $c' = 0$, whereas that of Storey (2003) requires that $c > 0$, $b = 0$, and $c' > 0$; he pointed out that a scientist can feasibly decide on $c/c'$ for a microarray experiment. If the net desirability is defined as $d(c, b, c') \equiv \sum_{i=1}^{m} d_i(c, b, c')$, then the Bayes error that Storey minimized, BE$(c, 0, c')$, is equal to $-E(d(c, 0, c'))/((c + c')m)$. Since maximizing $E(d(c, b, 0))$ and minimizing $-E(d(c, 0, b))$ yield the same optimal rejection region, the methods are mathematically equivalent. Thus, stating the problem in terms of false nondiscoveries as opposed to true discoveries is largely a matter of convention, but the latter is simpler and is more in line with the intuition that motivates reporting false discovery rates, which concerns false and true discoveries, not nondiscoveries. Unlike the technique used here, Storey's (2003) algorithm of optimization relies on a knowledge or estimation of probability densities, whereas Eqs. (20) and (19) only require the estimation of cumulative distributions. Scott and Berger (2003) and Müller *et al.* (2004) instead applied Bayesian methods to this optimization problem. An advantage of maximizing $\hat{D}(\tau)$ over these methods is that it does not require a Bayesian framework, but also applies to the case of fixed values of $H_i$.

Genovese and Wasserman (2001) considered a linear combination of the FDR and the *false nondiscovery rate* (FNR), and Storey (2003) optimized the corresponding linear combination the pFDR and the *positive false nondiscovery rate* (pFNR). The method of Storey (2003) would require fewer assumptions using the dFDR (12) and the *decisive false nondiscovery rate* (dFNR), the ratio of the expected number of false nondiscoveries to the expected number of nondiscoveries, or 0 if the expected number of nondiscoveries is 0. However, the maximization of $\hat{D}(\tau)$ or minimization of the estimated Bayes error is preferred since they are based on $c/b$, the ratio of the cost of a false discovery to the benefit of a true discovery or, equivalently, on $c/c'$. The linear combination methods, by contrast, are based on the ratio of the discovery rates. The problem is that the two ratios are not the same in general since the discovery rates have different denominators. Since it is difficult in most cases to specify and interpret a ratio of discovery rates, the straightforward approach of optimizing a total desirability or Bayes error is more widely applicable.

## 4.3 Test-dependent costs and benefits

The use of the dFDR and decision theory does not require all tests to have equal costs of false discoveries and benefits of true discoveries. If the costs and benefits associated with each test are equal to those of a sufficiently large number of other tests, as in the following application to microarray data, then the estimators of Eqs. (20), (19), and (22) and the maximization of $\hat{D}(\tau)$ can be applied separately to each subset of tests with the same costs and benefits, yielding a different rejection region for each subset. In the notation used above, $\forall_{k \in \{1,2,...,L\}} \forall_{i,j \in \iota_k} c_i = c_j, b_i = b_j,$ where $L \ll m$, $\forall_{j,k \in \{1,2,...,L\}} \iota_j \cap \iota_k = \emptyset,$ and $\bigcup_{k=1}^{L} \iota_k \equiv \{1, 2, ..., m\}$. It is assumed that each set $\iota_k$ has enough members that reliable estimates can be obtained by replacing $\{1, 2, ..., m\}$ with $\iota_k$ in the sums of Eqs. (20) and (19) and with test statistics of the null distribution generated from $\iota_k$. Let $\hat{D}_k(\tau_k)$ be the resulting estimate of the net desirability associated with $\iota_k$ for threshold $\tau_k$. Maximizing $\hat{D}_k(\tau_k)$ independently for each $\iota_k$ maximizes the estimated expected desirability, $\sum_{k=1}^{L} \hat{D}_k(\tau_k)$, but if a common rejection region is desired, then the estimate $\hat{D}_{b+c}(\tau)$ can be maximized instead by finding the optimal common threshold $\tau$, using $\hat{\pi}_0(b+c)$, a common estimate of $\pi_0$, found as described below. Since $\sum_{k=1}^{L} \hat{D}_k(\tau_k) \geq \hat{D}_{b+c}(\tau)$ if all estimates of $\pi_0$ and all null distributions are equal, the former approach tends to give a higher desirability, but the latter has an interesting connection with the weighted false discovery rate (WFDR) of Benjamini and Hochberg (1997). In analogy with the WFDR and the dFDR, define the *weighted decisive false discovery rate* (WdFDR) as

$$\Delta_w \equiv \begin{array}{l} \frac{E(\sum_{i=1}^{m} w_i V_i)}{E(\sum_{i=1}^{m} w_i R_i)}, \quad P(\sum_{i=1}^{m} w_i R_i > 0) \neq 0 \\ 0, \quad P(\sum_{i=1}^{m} w_i R_i > 0) = 0 \end{array}, \quad (25)$$

where $w$ is a vector of nonnegative weights, $(w_i)_{i=1}^{m}$. From Eq. (5), the expected net desirability can be expressed as

$$E(d) = b_1 E\left(\sum_{i=1}^{m} \left(\frac{b_i}{b_1} + \frac{c_i}{b_1}\right) R_i\right) \left(E\left(\sum_{i=1}^{m} \frac{b_i}{b_1} R_i\right) \bigg/ E\left(\sum_{i=1}^{m} \left(\frac{b_i}{b_1} + \frac{c_i}{b_1}\right) R_i\right) - \Delta_{b+c}\right), \quad (26)$$

Finding the threshold $\tau$ that maximizes $\hat{D}_{b+c}(\tau)$, the estimate of $E(d)$, thus yields an optimal estimate of $\Delta_{b+c}$:

$$\hat{\Delta}_{b+c}(\tau) \equiv \frac{\hat{\pi}_0(b+c) \sum_{i=1}^{M} (b_i + c_i) I(T_i^{(0)} \geq \tau) / M}{\sum_{i=1}^{m} (b_i + c_i) I(T_i \geq \tau) / m}, \tag{27}$$

with $\hat{\pi}_0(b+c)$ and $\hat{D}_{b+c}(\tau)$ similarly modified from Eqs. (19) and (22). Here, the null statistics are generated using a pool of all of the data, unlike the method in which each $\hat{D}_k(\tau_k)$ is independently maximized. Eqs. (26) and (27) indicate that $\hat{\Delta}_{b+c}(\tau)$ can be seen as an estimate of the WdFDR when each weight is set to the sum of the cost and benefit for the corresponding test.

In addition, there is a Bayesian interpretation of the WFDR. Given that the distribution of the test statistics is the mixture $F = \sum_{i=1}^{m} w_i F_j / \sum_{j=1}^{m} w_j$, the dFDR of test statistics distributed as $F$ [Eq. (16)] is equal to the WdFDR of test statistics distributed as $\sum_{i=1}^{m} F_j / m$:

$$\Delta(\tau) = \frac{\pi_0 (1 - F^{(0)}(\tau))}{1 - F(\tau)} = \frac{\pi_0 \left(1 - \sum_{i=1}^{m} w_i F_j^{(0)}(\tau) / \sum_{j=1}^{m} w_j\right)}{1 - \sum_{i=1}^{m} w_i F_j(\tau) / \sum_{j=1}^{m} w_j} =$$
$$\frac{\pi_0 \sum_{i=1}^{m} w_i (1 - F_j^{(0)}(\tau)) / \sum_{j=1}^{m} w_j}{\sum_{i=1}^{m} w_i (1 - F_j(\tau)) / \sum_{j=1}^{m} w_j} = \frac{E(\sum_{i=1}^{m} w_i V_i)}{E(\sum_{i=1}^{m} w_i R_i)} = \Delta_w \tag{28}$$

for $E(\sum_{i=1}^{m} w_i R_i) \neq 0$, and $\Delta(\tau) = \Delta_w = 0$ otherwise. Thus, the WdFDR is a posterior probability associated with $t \sim F$, according to Eq. (23). It follows that giving one test more weight than another is equivalent to increasing its representation in the sample by a proportional amount. Through the mixture $F$, the PFP, FDR, and pFDR have parallel relationships with their weighted counterparts.

Benjamini and Hochberg (1997) proposed the WFDR to take into account differing values or economic considerations among hypotheses, just as traditional Type I error rate control can give different weights for *p*-values or different significance levels to different hypotheses. However, without a compelling reason to have a common rejection region, optimiz-ing $\hat{D}_k(\tau_k)$ independently for each $\iota_k$ seems to better accomplish the goal of Benjamini and Hochberg (1997), in light of the above inequality.

## 5   Application to gene expression in cancer patients

The optimization of $\hat{\Delta}(\tau)$ was applied to the public microarray data set of Golub et al. (1999). Levels of gene expression were measured for 7129 genes in 38 patients with B-cell acute lymphoblastic leukemia (ALL) and in 25 patients with acute myeloid leukemia (AML). These preprocessing steps were applied to the raw "average difference" (AD) levels of expression: each AD value was divided by the median of all AD values of the subject to reduce subject-specific effects, then the normalized AD values were transformed by $(AD/|AD|)\ln(1+|AD|)$, as per Bickel (2002). The $i$th null hypothesis is that there is no systematic difference in transformed expression between groups in the $i$th gene; if it is rejected, then the gene is considered differentially expressed. The goal is to reject null hypotheses such that the expected desirability is maximized, given a cost-to-benefit ratio. That ratio is set to 19 to ensure that the probability of a false discovery of differential expression is no more than about 5%, according to equation (12). (Due to sampling error, there is no guarantee that the probabililty cannot exceed 5%.) The computations used the specific values $c = 19$ and $b = 1$, without loss of generality. (Alternately, one could estimate the cost of investigating a false discovery to the benefit of finding a truly differentially expressed gene.) To achieve the goal, the absolute value of a two-sample $t$-statistic,

$$T_i = \left| \mu_{\text{ALL},i} - \mu_{\text{AML},i} \right| / (\sigma^2_{\text{ALL},i}/38 - \sigma^2_{\text{ALL},i}/25)^{1/2}, \tag{29}$$

was computed for each gene for both the original preprocessed data and for $M = 1000$ random permutations of the columns of the $7129 \times (38 + 25)$ data matrix, following the unbalanced permutation method used by Dudoit et al. (2002). (The whole columns were permuted, instead of the expression values for each gene independently, to preserve the gene-gene correlations.) Finally, the threshold $\tau$ was found such that $\hat{D}(\tau)$ is maximized, using the estimates computed from Eqs. (20), (19), and (22). Fig. 2 displays $\hat{D}(\tau)$ for this data. Table 2 compares the results to those obtained by choosing $\tau$ such that the highest number of discoveries were made with the constraint $\hat{\Delta}(\tau) \leq 5\%$ as dFDR control. Since the method of Benjamini and Hochberg (1995) effectively uses $\pi_0 = 1$ (Storey 2002a), as does the algorithm of Bickel (2002, 2004a), Table 2 also shows the effect of conservatively replacing $\hat{\pi}_0$ with 1.

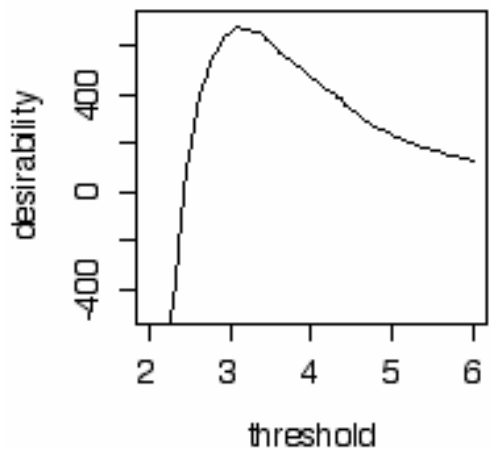

Estimated desirability $\hat{D}(\tau)$ for the ALL-AML comparison (22), using the estimate $\hat{\pi}_0 \approx 0.59$, from Eq. (19).

**Figure 2**

|  | $\max \hat{D}(\tau)$ $P(H_i = 0 \mid t_i) \lesssim 5\%$ $\hat{\pi}_0 \approx 0.59$ | $\max \hat{D}(\tau)$ $P(H_i = 0 \mid t_i) \lesssim 5\%$ $\hat{\pi}_0 \equiv 1$ | $\hat{\Delta}(\tau) \leq 5\%$ $\hat{\pi}_0 \approx 0.59$ | $\hat{\Delta}(\tau) \leq 5\%$ $\hat{\pi}_0 \equiv 1$ |
|---|---|---|---|---|
| Threshold, $\tau$ | 3.14 | 3.37 | 2.44 | 2.73 |
| #discoveries | 910 | 768 | 1496 | 1212 |
| Desir., $\hat{D}(\tau)$ with $\hat{\pi}_0 \approx 0.59$ | **683** | 656 | 1.4 | 500 |
| Desir., $\hat{D}(\tau)$ with $\hat{\pi}_0 \equiv 1$ | 524 | **578** | −1043 | 2.3 |
| dFDR, $\hat{\Delta}(\tau)$ with $\hat{\pi}_0 \approx 0.59$ | 0.0125 | 0.0073 | **0.0500** | 0.0294 |
| dFDR, $\hat{\Delta}(\tau)$ with $\hat{\pi}_0 \equiv 1$ | 0.0212 | 0.0124 | 0.0849 | **0.0499** |

Comparison of estimates of the proposed method of the first column (maximizing $\hat{D}(\tau)$ with $\hat{\pi}_0$ from Eq. (19)) to three alternate methods. Each number in **bold face** is a maximal value constrained according to the method of its column.

**Table 2**

Comparisons between other groups can easily be added to these results. For example, one might be interested not only in knowing which genes are differentially expressed between B-cell ALL and AML patients, but also which ones are differentially expressed between the B-cell ALL group and the 9 patients of Golub et al. (1999) with T-cell ALL. A family-wise error rate control approach, such as that used by Dudoit et al. (2002), would require additional adjustments to the *p*-values, with the result that fewer genes would be found to be differentially expressed as the number of group-group comparisons increases, but the approach of ingent rejection regions, given by $\tau_1$ and $\tau_2$, does not require any modification to the results of Table 2. Suppose that the benefit of a true discovery of differential expression between the B-cell ALL and T-cell ALL groups is twice that between the B-cell ALL and AML groups, and that the cost of investigating a false discovery is the same for all comparisons. Then the only difference in the method applied to the new set of comparisons is that the benefit of each true discovery is set to 2 instead of 1. In the above notation, $\iota_1$ specifies the *t*-statistics corresponding to the comparisons to the AML group and $\iota_2$ specifies those corresponding to the comparisons to the T-cell ALL group. In microarray experiments, this kind of multiplicity of subject groups as well as a multiplicity of genes is common, and the dFDR and decision theory together provide a flexible, powerful framework for such situations. The computed threshold is $\tau_2 \approx 3.64$, yielding 350 genes considered differentially expressed and a dFDR estimate of 0.0418; these values correspond to $\tau_1 \approx 3.14$, 910 genes, and a dFDR estimate of 0.0125 from the first column of Table 2. The two thresholds are close enough that forcing a common threshold, as in Eqs. (27) and (28), would not have a drastic effect on the conclusions drawn in this case. By contrast, the marked difference between the two dFDR estimates is a reminder that traditional FDR control does not maximize the desirability. Even so, each estimate is controlled in the sense of satisfying (12).

In estimating the pFDR, Storey (2003) considered the expression patterns of genes to be weakly dependent, but that may presently be impossible to test statistically since the number is genes is typically much greater than the number of cases. Storey (2002b) hypothesized approximate weak dependence since genes work in biochemical pathways of small groups of interacting genes. However, such weak dependence would only follow if there were no interaction between one pathway and another, an assumption that often does not hold. Indeed, the expression of thousands of genes in eukaryotes is controlled by highly connected networks of transcriptional regulatory proteins (Lee et al., 2002). Thus, we hypothesize instead that while weak dependence might be a good approximation for prokaryotes, there is stronger dependence between eukaryotic genes. Fractal (scale-free, self-similar) statistics concisely describe biological long-range correlations (Bickel and West, 1998) and may aid the study of genome-wide interactions as more data become available. Agreeably, gene networks reconstructed from expression data tend to have fractal connectivity distributions (Agrawal, 2002; Bickel, 2004b). A further challenge is that having only a few microarrays per group can cause marked violations of weak dependence, for example, when any array effects not removed by normalization or when any other substantial random effects are present. Mixed-effects models may yield p-values without such problems, but such models lack statistical power for small numbers of microarrays.

However, the inability to make assumptions about the dependence of hypothesis tests need not prevent the reliable use of false discovery rates. Benjamini and Yekutieli (2001) proposed a method of controlling the false discovery rate without such assumptions, but it has a lower rate of true discoveries. Another appoach to the problem of general dependence uses conventional methods to control or estimate a redefined false discovery rate (Bickel 2004d).

## Acknowledgements

I thank Dan Nettleton for bringing previous work on the PFP to my attention, and I thank the two anonymous reviewers for helpful comments.